\documentclass[final,onefignum,onetabnum]{siamart171218}

\usepackage[paperwidth=7in, paperheight=10in, margin=.875in]{geometry}
\usepackage{amsfonts,amssymb}
\usepackage{amsmath}
\usepackage{graphicx}
\usepackage{booktabs}
\usepackage{cite}
\usepackage{enumerate}
\renewcommand{\theequation}{\arabic{section}.\arabic{equation}}

\usepackage{amsmath,amssymb,cleveref,graphicx,subcaption,amsbsy}
\usepackage[disable]{todonotes}

\usepackage[numbers]{natbib}

\DeclareMathOperator*{\argmin}{argmin}

\newcommand{\ext}[1]{\underset{#1}{\text{ext}}}
\graphicspath{{/figures/}}

\newcommand{\B}[1]{{\boldsymbol #1}}
\newcommand{\bx}{\boldsymbol x}

\newcommand{\R}{\mathbb{R}}

\newcommand*{\defeq}{\mathrel{\vcenter{\baselineskip0.5ex \lineskiplimit0pt
                     \hbox{\scriptsize.}\hbox{\scriptsize.}}}%
     		     =}

\def\XXint#1#2#3{{\setbox0=\hbox{$#1{#2#3}{\int}$ }
\vcenter{\hbox{$#2#3$ }}\kern-.6\wd0}}

\title{A Model for Optimal Human Navigation with Stochastic Effects\thanks{Submitted: \today}}

\author{Christian Parkinson\thanks{Department of Mathematics, UCLA, Los Angeles, CA 90095, (chparkin@math.ucla.edu).}
\and David Arnold\thanks{Department of Mathematics, UCLA, Los Angeles, CA 90095, (darnold@math.ucla.edu). Department of Mathematics and Statistics, Macquarie University, Sydney 2109, Australia (david.arnold@mq.edu.au).}
\and Andrea L. Bertozzi\thanks{Department of Mathematics, UCLA, Los Angeles, CA 90095, (bertozzi@math.ucla.edu).}
\and Stanley Osher\thanks{Department of Mathematics, UCLA, Los Angeles, CA 90095, (sjo@math.ucla.edu).}}

   \allowdisplaybreaks

\begin{document}

\maketitle

\begin{abstract}
We present a method for optimal path planning of human walking paths in mountainous terrain, using a control theoretic formulation and a Hamilton-Jacobi-Bellman equation. Previous models for human navigation were entirely deterministic, assuming perfect knowledge of the ambient elevation data and human walking velocity as a function of local slope of the terrain. Our model includes a stochastic component which can account for uncertainty in the problem, and thus includes a Hamilton-Jacobi-Bellman equation with viscosity. We discuss the model in the presence and absence of stochastic effects, and suggest numerical methods for simulating the model. We discuss two different notions of an optimal path when there is uncertainty in the problem. Finally, we compare the optimal paths suggested by the model at different levels of uncertainty, and observe that as the size of the uncertainty tends to zero (and thus the viscosity in the equation tends to zero), the optimal path tends toward the deterministic optimal path. 
\end{abstract}

 


\section{Introduction} \label{intro} Optimal path planning is a classical problem in control theory and engineering. The basic mathematical formulation of the problem includes an agent whose movement is constrained by an equation of motion, and who desires to optimally travel from one point to another while obeying this equation. We will consider optimal path planning in the context of hikers traversing mountainous terrain. 

Early approaches to the optimal path planning problem date back to the 1950's and were originally discrete in nature: discretizing the continuous domain into a weighted graph and then employing Dijkstra's algorithm and its many variants \cite{Dijkstra}. Several of these discrete or semi-discrete approaches are still being developed and improved today \cite{Hirsch, JSBMitchell, Papadimitrou,krishnaswamy1995resolution, SARANYA2016178, RealTime}. However, since the late 1990's, there has also been significant interest in solving the problem using variational and partial different equation (PDE) based models. Sethian and Vladimirsky determine optimal paths on manifolds by formulating the dynamic programming principle and Hamilton-Jacobi-Bellman (HJB) equation as a continuous version of Dijkstra's algorithm \cite{SethVlad2,SethVlad1}. In a similar approach, one can use the level set method of Osher and Sethian \cite{OsherSethian} to resolve optimal paths \cite{CecilLevelSet, OptPathPaper}. Alternatively, Martin and Tsai suggest a steady-state HJB formulation for determining optimal paths on manifolds which are represented computationally by unstructured point clouds \cite{MartinTsai}. In application, variational methods for optimal path planning are extensively used in the context of self-driving vehicles. This problem was first considered by Dubins \cite{Dubins} in discrete form, but later reformulated continuously, and adapted to answer a number of questions involving impassable obstacles and reachable sets, among others \cite{AgarwalWang,Lygeros,Takei,TomlinAircraft}. Human movement has also been modeled using HJB equations, whether it be walking while expending minimal energy \cite{BipedWalking,Papadopoulos2016}, reach-avoid games like capture-the-flag \cite{TomlinReachAvoid, SurvEvas}, or pedestrian flow modeling \cite{AnisPedestrian}.

One feature of all the models referenced above is that they are completely deterministic. For example, in the simplest reach-avoid games, the strategy of each team is known to the opposing team. Assumptions like this may not be realistic in practice, and thus it is important to incorporate some uncertainty into the models, and this can be done in a number of ways. In the case of reach-avoid games, Gilles and Vladimirsky suggest paths for the attackers or defenders that minimize or maximize risk in different ways \cite{SurvEvasUncertain}. In optimal path planning, Shen and Vladimirsky account for weather effects on a sailboat by designing a piecewise deterministic algorithm, wherein the travel velocity changes randomly but at discrete times \cite{ShenVladimirsky}. These type of stochastic effects are of special interest to those working on self-driving vehicles wherein misaligned axles, miscalibrated sensors or a host of other variables could significantly perturb an optimal driving path, and these can be modeled using randomness in a number of different ways \cite{Aurnhammer2001, LimRus, Luna, Marti, Marti2015}. Others have suggested optimal path planning for underwater unmanned vehicles using a stochastic drift term to account for the effects of the ocean's current \cite{Lolla,Subramani1,Subramani2}.

In this paper, we propose a method for optimal path planning where the walking velocity depends on the local slope in the direction of motion (and is thus anisotropic), and there is uncertainty in the equation of motion of the traveler. That is, rather than an ordinary differential equation, our equation of motion is a stochastic differential equation where Brownian motion can account for uncertainty in human walking speed, the ambient elevation data, terrain traversibility or other uncertainties. In doing so, we follow a similar formulation as in a previous paper \cite{OptPathPaper}. However, since the Hamilton-Jacobi-Bellman equation for the stochastic case has viscosity, the level set method is no longer applicable, so we opt for a model more rooted in control theory. This paper is organized as follows: in section \ref{model}, we discuss the setup of our model, including some of the underlying mathematical formalism and how it pertains to optimal path planning. In section \ref{numerics}, we discuss the numerical methods which we used to simulate the model. In section \ref{results}, we discuss the results of our simulations. Specifically, we test our model against both synthetic data and real elevation data taken from the area surrounding El Capitan, a mountain in Yosemite National Park. We compare two different notions of optimal paths when stochastic effects are present. Finally, we observe that as the size of the uncertainty tends to zero (and thus the solution of the stochastic HJB equation tends to the viscosity solution of the ordinary HJB equation \cite{CrandallLions}), the optimal path tends toward the deterministic optimal path.

\section{Mathematical Model}\label{model} We discuss the construction of our model, only briefly mentioning the underlying control theoretic concepts. Treatments of these concepts with varying levels of rigor can be found in several books \cite{Bertsekas,Carmona,Fleming,Pham,Nisio}. The main insights here are how the controlled equation of motion leads to a Hamilton-Jacobi-Bellman equation, and how we can use this equation to determine the optimal control. 

\subsection{The Deterministic Optimal Path Planning Model}\label{HJBDeterministic} We consider the problem of planning optimal walking paths in mountainous terrain. We imagine that a hiker is standing at a point $x_{0} \in \R^2,$ and wishes to walk to a point $x_{\text{end}} \in \R^2$. We are given the elevation profile $E: \R^2 \to \R$ which provides the elevation at any point in space, and velocity function $V: \R \to [0,\infty)$ which measures walking velocity as a function of slope. Here slope is in units of grade, so that slope of $1$ corresponds to a $45^\circ$ incline, and slope of $-1$ corresponds to a $45^\circ$ decline. Let $\B x: [0,T] \to \R^2$ represent the current position of the hiker, where $T > 0$ is the total walking time. Note, this terminal time $T$ is a parameter which one specifies at the outset. Our control variable will be the walking direction along the path, $\B s: [0,T] \to \mathbb S^1 \defeq \{a \in \mathbb R^2 \, : \, \| a\|_2 = 1\}$. Given all of this, the equation of motion for our hiker is \begin{equation} \label{eq:contODE} \begin{split} &\dot{\B x}(t) = V(\nabla E(\B x(t)) \cdot \B s(t)) \B s(t), \,\,\,\,\, 0 < t \le T, \\ 
&\B x(0) = x_0.
\end{split}\end{equation} Note that $\nabla E(\B x) \cdot \B s$ represents the local slope in the walking direction. 

For the cost functional, we simply consider Euclidean distance to the end point. Given a particular control parameter $\B s:[0,T] \to \mathbb S^1$, define
\begin{equation} \label{eq:cost} C[\B s(\cdot)] = \| \B x(T) - x_{\text{end}}\|_2, \end{equation} where $\B x$ is constrained by \eqref{eq:contODE}. To arrive at the Hamilton-Jacobi-Bellman equation, we invoke the dynamic programming principle, stating that to find a globally optimal path, we need only consider paths which are locally optimal. That is, for any given $x \in \R^2$ and $t \in (0,T]$, one can consider optimal paths on the interval $(t,T]$ with $\B x(t) = x$. We define the cost function by \begin{equation}\label{eq:costFunc} u(x,t) = \inf_{\B s } \big\{C_{x,t}[\B s(\cdot)] \big\}, \end{equation} where $C_{x,t}[\B s(\cdot)]$ denotes the same cost functional, but restricted to paths $\B x$ on the interval $(t,T]$ with $\B x(t) = x$. Provided that the map $x \mapsto V(\nabla E(x) \cdot s)$ is sufficiently regular \cite{Fleming}, one can show that the cost function $u$ satisfies the equation \begin{equation} \label{eq:HJB} \begin{split}  
&u_t(x,t) + \inf_{s \in \mathbb S^1} \left\{ V(\nabla E(x) \cdot s) [\nabla u(x,t) \cdot s ]\right\} = 0, \\
&u(x,T) = \| x - x_{\text{end}}\|_2.
\end{split}\end{equation} This is the Hamilton-Jacobi-Bellman (HJB) equation for this control problem. Note that this is a terminal value problem, for which we have data at time $T$ and we solve backwards in time on the interval $[0,T)$. This correspondence between the control problem and the HJB equation will hold under the condition that $x \mapsto V(\nabla E(x) \cdot s)$ is Lipschitz with a constant which is uniform over $s \in \mathbb S^1$ \cite{Fleming}. In our application, this will likely hold, but we cannot guarantee it because we will use real elevation data $E(x)$. Thus this may only be a formal correspondence, though empirically, the method appears to work even when the elevation data is somewhat non-smooth. Another note is that the infimum in \eqref{eq:HJB} is being taken over a compact set, so again, under some mild assumptions of continuity, this will be realized as a minimum.

The question remains of how to use \eqref{eq:HJB} to determine the optimal control, and thus the optimal trajectory. At each point, the optimal control is given by the infimum in \eqref{eq:HJB}. Thus if one can solve \eqref{eq:HJB}, the optimal control $\B s^*$ is \begin{equation} \label{eq:controlMin}\B s^*(x,t) = \argmin_{s \in \mathbb S^1} \left\{ V(\nabla E(x) \cdot s) [\nabla u(x,t) \cdot s ]\right\}.\end{equation}  When this minimization problem is solved, the optimal trajectory is then resolved by solving the differential equation \begin{equation} \label{eq:actualODE}\begin{split} &\dot{\B x}(t) = V(\nabla E(\B x(t)) \cdot \B s^*(\B x(t),t)) \B s^*(\B x(t),t), \,\,\,\,\, 0 < t \le T, \\ 
&\B x(0) = x_0. \end{split} \end{equation}

It was mentioned above that the terminal time $T$ is a parameter which one must specify beforehand. It is important to note what this means for the model. Notice that for our cost functional, we have chosen a lump sum cost which is the distance from the end of the path $\B x(T)$ to the desired end point $x_{\text{end}}$. The optimal control problem attempts to minimize this cost. Thus the path that we observe will be the path $\B x$ whose end point is as close to $x_{\text{end}}$ as possible, given time $T>0$. If we select $T$ too small, the path will not reach the end point, and in some cases this can lead to some interesting decisions regarding how the path should be constructed. We discuss this further in section \ref{theRoleOfT}.

\subsection{The Stochastic Optimal Path Planning Model} \label{HJBStochastic} In the above, we assume that all data is known perfectly and that there is no uncertainty. In reality, weather effects, instrumentation noise, incomplete elevation data or any number of other things could cause uncertainty in the equation of motion. Accordingly, we can account for random effects by considering a stochastic equation of motion \begin{equation} \label{eq:SDE}\begin{split} &d\B X_{t} = V(\nabla E(\B X_t) \cdot \B s(t)) \B s(t)dt + \sigma(\B X_t, \B s(t)) d \B W_t, \,\,\,\,\,\,\, 0 < t \le T \\ & \B X_0 = x_0, \end{split} \end{equation} where $\B W_t$ is an ordinary two-dimensional Brownian motion whose coordinates $W^{(1)}_t$ and $W^{(2)}_t$ are independent one-dimensional Brownian motions, and $\sigma$ is some function that determines the uncertainty. 

In this case, we define the expected cost function \begin{equation} \label{eq:costStoch} C[\B s(\cdot)] = \mathbb E \left( \| \B X_T - x_{\text{end}}\|_2 \right),\end{equation} and perform similar steps as above to arrive at a stochastic version of the Hamilton-Jacobi-Bellman equation. That is, we localize the problem and define \begin{equation}\label{eq:costFuncStoch} u(x,t) = \inf_{\B s }  \left\{ \mathbb E \big(C_{x,t}[\B s(\cdot)] \big) \right\}, \end{equation} where again  $C_{x,t}[\B s(\cdot)]$ denotes the cost functional restricted to trajectories $\B X$ on the interval $(t,T]$ with $\B X_t = x$.  A key point in the derivation of the HJB equation invokes the chain rule for the function $t \mapsto u(\B x(t),t)$, which explains the appearance of $u_t$ and $\nabla u$ in the HJB equation. We are now concerned with the map $t\to u(\B X_t, t)$ and we must consider second order derivatives because of the fundamental relationship $(dW^{(i)}_t)^2 \sim dt$ for $i =1,2$. Thus the stochastic Hamilton-Jacobi-Bellman equation reads \begin{equation} \label{eq:SHJB} \begin{split}  
&u_t(x,t) + \inf_{s \in \mathbb S^1} \left\{ V(\nabla E(x) \cdot s) [\nabla u(x,t) \cdot s ] + \frac{1}{2} \sigma(x,s)^2 \Delta u(x,t) \right\} = 0, \\
&u(x,T) = \| x - x_{\text{end}}\|_2.
\end{split}\end{equation} Again, this is a terminal value problem, solved backwards from $t = T$ to $t =0$. Thus the positive sign on the diffusion is the correct sign to ensure that the diffusion has a smoothing effect on the solution, and there is no danger of finite-time blow-up.  

As above, to determine the optimal control, we solve \eqref{eq:SHJB}, and at each point set \begin{equation} \label{eq:controlMinStoch}\B s^*(x,t) = \argmin_{s \in \mathbb S^1} \left\{ V(\nabla E(x) \cdot s) [\nabla u(x,t) \cdot s  + \frac{1}{2} \sigma(x,s)^2 \Delta u(x,t)]\right\}.\end{equation} At this point, there is a choice as to how to construct the path. One can use the optimal control computed from the stochastic HJB equation but simulate the deterministic path equation \begin{equation} \label{eq:actualODEStoch}\begin{split} &\dot{\B x}(t) = V(\nabla E(\B x(t)) \cdot \B s^*(\B x(t),t)) \B s^*(\B x(t),t), \,\,\,\,\, 0 < t \le T, \\ 
&\B x(0) = x_0, \end{split} \end{equation} to arrive at an optimal path. Alternatively, to compute a single instance of a path, we could simulate the equation \begin{equation} \label{eq:actualSDE}\begin{split} &d\B X_t = V(\nabla E(\B X_t) \cdot \B s^*(X_t, t)) \B s^*(X_t, t)dt + \sigma(\B X_t, \B s^*(\B X_t,t)) d \B W_t, \,\,\,\,\,\,\, 0 < t \le T \\ & \B X_0 = x_0. \end{split} \end{equation} Because this is only one instance and is subject to randomness, the path will not necessarily connect the points $x_0$ and $x_{\text{end}}$. However, we can average over many realizations to arrive at an expected optimal path. This is summarized in table~\ref{fig:twoMethods}.

\begin{table}[htbp]
\centering
\caption{The two options for how to compute optimal paths with uncertainty.}
\begin{tabular}{l l}\toprule
Method $(i)$ & Method $(ii)$ \\
\midrule
Stochastic HJB Equation \eqref{eq:SHJB}& Stochastic HJB Equation \eqref{eq:SHJB} \\
\hspace{1.9cm}$\to$ control values & \hspace{1.9cm}$\to$ control values\\ \midrule
Deterministic equation of motion \eqref{eq:actualODEStoch}  & Stochastic equation of motion \eqref{eq:actualSDE}\\ 
 \hspace{1cm}$\dot{\B x} = V(\nabla E(\B x) \cdot \B s) \B s$ & $\mathbb E \big( \,\,\,\, d\B X_t = V(\nabla E(\B X_t) \cdot \B s) \B s dt + \sigma d \B W_t\,\,\,\, \big)$\\
\hspace{1.9cm}$\to$ optimal path & \hspace{1.9cm}$\to$ expected optimal path \\
\bottomrule
\end{tabular}
\label{fig:twoMethods}
\end{table}

These methods have different interpretations and may be better suited to modeling different physical scenarios. In using method $(i)$, the uncertainty is factored into the planning of the route, but upon traversing the route, there is no uncertainty in the velocity. This could model a hiker walking through a forest. The hiker does not feel random perturbations in the walking velocity at each step; rather, the uncertainty is in the exact form of the landscape that lies ahead. Method $(ii)$ may be of more practical use to a company shipping goods from one port to another, wherein each boat that makes the trip will actually feel random perturbations in velocity due to wind or currents. We will use both methods to compute paths and compare the results in section~\ref{results}.

\subsection{Our Model} \label{ourmodel} In order to simulate the model, we simply need to determine some parameter values, specifically those of $E,V$ and $\sigma$. For the elevation data $E$, we will begin by using synthetic elevation data which is mostly flat but with a few ``mountains" included which we would expect the hiker to avoid. After this, we will use real elevation data in the area surrounding El Capitan, a large granite cliff face in Yosemite National Park in California. We will specify the elevation data which we use for each simulation in section \ref{results}.

For the walking velocity function, we use a modified version of the function of Irmischer and Clarke \cite{IC2017}. Our specific velocity function is \cite{OptPathPaper}: \begin{equation}\label{eq:VelUs} V(S) = 1.11 \exp\left(-\frac{(100S + 2)^2}{2345}\right). \end{equation} Additionally, if we use this velocity function and only consider the walking direction $\nabla E(x)\cdot s$, then the walking velocity in the tangential and normal directions to the path are completely decoupled. Thus one could walk easily in the east-west direction, even when the grade in the north-south direction is arbitrarily steep. To prevent this, we penalize the velocity if the grade is sufficiently steep in any direction \cite{OptPathPaper}.

Lastly, one must determine the exact form of the uncertainty $\sigma$ in the stochastic equation of motion. For our purposes, we take $\sigma$ constant, so that we model uncertainty in the walking velocity in a general sense without specifying the exact nature of the uncertainty. As a consequence, if we reconsider \eqref{eq:SHJB}, we notice that the viscosity term is independent of the control variable $s$, and thus the equation can be re-written \begin{equation} \label{eq:SHJBOurs} \begin{split}  
&u_t(x,t) + \frac{\sigma^2}{2}\Delta u(x,t) + \inf_{s \in \mathbb S^1} \left\{ V(\nabla E(x) \cdot s) [\nabla u(x,t) \cdot s ] \right\} = 0, \\
&u(x,T) = \| x - x_{\text{end}}\|_2.
\end{split}\end{equation} This case is interesting because now the optimal control is resolved exactly as in the deterministic case, and \eqref{eq:SHJBOurs} is reminiscent of the viscous Hamilton-Jacobi equation considered by Crandall and Lions \cite{CrandallLions} when establishing the vanishing viscosity method for Hamilton-Jacobi equations. Thus if our Hamiltonian \begin{equation} \label{eq:Hamiltonian} H(x,p) \defeq \inf_{s \in \mathbb S^1} \left\{ V(\nabla E(x) \cdot s) [ p \cdot s ] \right\}\end{equation} is continuous, we expect that the solution $u^{(\sigma)}$ to equation \eqref{eq:SHJBOurs} will converge to the viscosity solution $u$ of \eqref{eq:HJB} as $\sigma \searrow 0$ \cite{CrandallLions}. Again, continuity of $H(x,p)$ will depend on the nature of the elevation data, but we can observe this convergence empirically by considering the optimal path constructed by the the method at varying levels of uncertainty $\sigma$. 

Other modeling decisions could be made. For example one could consider a running cost, allow for infinite horizon time, or allow $\sigma$ to depend on $x$ and $s$ (perhaps accounting for local slope as with the velocity function). We chose our cost functional and finite horizon time in analogy with the previous level set model for deterministic path planning \cite{OptPathPaper}.  With these decisions, the HJB equation remains time-dependent, and can thus be approximated very simply at high-order accuracy using techniques like ENO and WENO \cite{WENO,OsherShu}. As stated above, we will discuss the ramifications of this decision in section \ref{theRoleOfT}. Likewise, choosing $\sigma$ constant simplifies the numerics in that it allows the diffusion to be resolved implicitly as seen in section \ref{sec:semiImp}. Alternatively, if one wishes to have $\sigma(x,s)$, it will likely be difficult to resolve the diffusion implicitly, and a steady-state, infinite horizon time formulation like those in \cite{MartinTsai, TomlinReachAvoid} may be necessary. This would also require more involved numerical schemes; one could likely use a modified fast marching \cite{Tsitsiklis,SethVlad2,SethVlad1} or fast sweeping method \cite{kao2004lax,TsaiOsherSweep}. Indeed, in the latter reference \cite{TsaiOsherSweep}, a fast sweeping scheme is specifically applied to an anisotropic Eikonal equation which describes geodesic distance on the graph of a smooth function.

\section{Numerical Methods}\label{numerics} There are several numerical concerns that must be addressed in order to simulate these equation; both general concerns for numerical Hamilton-Jacobi equations and specific concerns relating to our model. We begin by discussing general notions for solving Hamilton-Jacobi type equations numerically.

\subsection{An Explicit Scheme for \eqref{eq:HJB}} We consider a general Hamilton-Jacobi equation of the form \begin{equation}\label{eq:HJ}
u_t + H(u_x,u_y) = 0, \,\,\,\,\, x \in \mathbb R^2, \,\, t > 0, \end{equation} with initial data at $t=0$. Our Hamiltionian \eqref{eq:Hamiltonian} depends also on $x$, but for the sake of numerical methods, this is unimportant, so we suppress it here to simplify notation. Since solutions of Hamilton-Jacobi equations have kinks \cite{Bardi1997,CrandallLions,osher2003level}, na\"ive forward, backward and centered differencing methods may fail to accurately simulate the equation. Accordingly, in a numerical scheme, one must replace the Hamiltonian $H(u_x, u_y)$ with a \emph{numerical Hamiltonian} $\hat H(u_x^+,u_x^-,u_y^+,u_y^-)$ which deftly combines the forward differences $u_x^+, u_y^+$ and backwards differences $u_x^-, u_y^-$ so as to smooth the solution, minimize oscillation near kinks, track the characteristics, or otherwise capture the dynamics of the equation. 

Following Osher and Shu \cite{OsherShu}, we use the Godunov numerical Hamiltonian \begin{equation}\hat H_G(u_x^+,u_x^-,u_y^+,u_y^-) = \ext{u \in I(u_x^-,u_x^+)}\,\,\, \ext{v \in I(u_y^+,u_y^-)} H(u,v) \end{equation} where \begin{equation}I(a,b) = [\min(a,b), \max(a,b)] \hspace{0.5cm} \text{ and } \hspace{0.5 cm}\ext{x \in I(a,b)} = \left \{\begin{matrix} \min_{a \le x \le b} & \text{if } a \le b, \\ \max_{b \le x \le a} & \text{if } a > b. \end{matrix} \right. \end{equation} These minima and maxima are designed to track the characteristics of the equation, thus accurately approximating the Hamiltonian without resorting to excessive numerical diffusion as is present in the Lax-Friedrichs scheme; we comment on this further later. Osher and Shu \cite{OsherShu} suggest a method for approximating $u_x$ and $u_y$ with forward and backward difference schemes which are accurate to arbitrarily high order, though we opt for  first order approximations which suffice in this application. In the absensce of viscosity, we then explicitly integrate \eqref{eq:HJ} in time using explicit Euler time stepping. 

Specifically while \eqref{eq:HJB} has infinite spatial domain, for the purpose of simulating it, we draw some box $[a,b] \times [c,d]$ containing the points $x_0$ and $x_{\text{end}}$. We also invert the time variable, so that we instead solve for $u(x,y,T-t)$ on the interval $t \in (0,T]$. Next discritize this box, and the time interval uniformly: \begin{equation} \label{eq:disc}\begin{split} x_i &= a + i \frac{(b-a)}{N},   \,\,\,\,\,\, i = 0,1,\ldots, N, \\ y_j &= c + j \frac{(d-c)}{M}, \hspace{1.5mm} j = 0,1,\ldots,M, \\ t_k &= \frac{kT}{K}, \hspace{1.3cm} k = 0,1,\ldots,K. \end{split} \end{equation} Then if $u^k_{ij} \defeq u(x_i,y_j,T-t_k)$, our explicit time stepping scheme is \begin{equation} \label{eq:explicitScheme} u_{ij}^{k} = u^{k-1}_{ij} + \Delta t \hat H_G(u_x^+,u_x^-,u_y^+,u_y^-)^{k-1}_{ij}\end{equation} for $i=1,\ldots, N-1, j= 1,\ldots, M-1$ and $k=1,\ldots,K$, where $\Delta t \defeq T/K$. Note, the positive sign in front of $\hat H_G$ in \eqref{eq:explicitScheme} is due to the time inversion; this would ordinarily be negative when moved to the right hand side of \eqref{eq:HJB}. If one is resolving the spatial derivatives to higher order accuracy and desires to maintain this accuracy, one can easily replace the explicit Euler time integration with a higher order Runge-Kutta scheme \cite{OsherShu}.

The layers of nodes corresponding to $i=0, i=N, j=0$ or $j=M$ are necessarily boundary layers, because for example, we cannot compute the backwards difference approximation $u_x^-$ at the nodes where $i=0$. Thus after evaluating \eqref{eq:explicitScheme} at each time step, we must enforce some artificial boundary condition, which is not prescribed in the differential equation, but is rather a purely numerical consideration. We use the boundary conditions suggested by Kao et al. \cite{kao2004lax} which extrapolate while also attempting to minimize oscillation and prevent information from entering through the boundary: \begin{equation}\label{eq:BCs} \begin{split}
u_{0,j}^k &= \min(\max(2u^k_{1,j}-u^k_{2,j},u^k_{2,j}),u^{k-1}_{0,j}),\\ 
u_{N,j}^k &= \min(\max(2u^k_{N-1,j}-u^k_{N-2,j},u^k_{N-2,j}),u^{k-1}_{N,j}),\\
u_{i,0}^k &= \min(\max(2u^k_{i,1}-u^k_{i,2},u^k_{i,2}),u^{k-1}_{i,0}),\\
u_{i,M}^k &= \min(\max(2u^k_{i,M-1}-u^k_{i,M-2},u^k_{i,M-2}),u^{k-1}_{i,M}).
\end{split} \end{equation} 

For our particular application, the maximum velocity at which information flows along characteristics is $V_{\text{max}} \defeq \max_S V(S) = 1.11$, and so this scheme will be stable so long as the parameters $(\Delta x, \Delta y, \Delta t) \defeq \left(\tfrac{b-a}{N},\tfrac{d-c}{M}, \tfrac T K\right)$ satisfy the CFL condition \cite{osher2003level}: \begin{equation} \label{eq:CFL} \Delta t \cdot V_{\text{max}}\left( \frac{1}{\Delta x} + \frac{1}{\Delta y}\right) < 1. \end{equation} 

\subsection{A Semi-Implicit Scheme for \eqref{eq:SHJBOurs}} \label{sec:semiImp}In order to numerically simulate the reaction-diffusion equation \eqref{eq:SHJBOurs}, one could simply insert the centered difference approximation to $\Delta u$ and continue with the explicit Euler time stepping. However, this will require exceedingly small time discretization, since the stability condition for forward Euler time stepping for a diffusion operator is of the form $\Delta t = O((\Delta x)^2, (\Delta y)^2)$. Instead, we resolving the diffusion implicitly. Thus, our scheme for \eqref{eq:SHJBOurs} is \begin{equation} \label{eq:semiImplicitScheme} u_{ij}^{k} - \frac{\sigma^2 \Delta t}{2\Delta x} (u^+_x - u^-_x)^k_{ij} -  \frac{\sigma^2 \Delta t}{2\Delta y} (u^+_y - u^-_y)^k_{ij} = u^{k-1}_{ij} + \Delta t \hat H_G(u_x^+,u_x^-,u_y^+,u_y^-)^{k-1}_{ij}.\end{equation} 

Since implicit Euler time stepping for diffusion is unconditionally stable, our discretization is still only bound by the CFL condition \eqref{eq:CFL}. For larger values of $\sigma$, the resulting diffusion will smooth the solution $u$, and thus sophisticated numerical Hamiltonians are probably no longer necessary. However, we have implemented this scheme as stated, so that as $\sigma \searrow 0$, the semi-implicit scheme \eqref{eq:semiImplicitScheme} reverts to the explicit scheme \eqref{eq:explicitScheme}, and there are no stability issues. We note that this semi-implicit scheme is only available when the uncertainty $\sigma$ in the equation of motion is independent of the control variable. If $\sigma$ depends on $s$, then the Hamilton-Jacobi-Bellman equation takes the form \eqref{eq:SHJB}, and in this case the Hamiltonian cannot be decoupled from the diffusion term. Thus one would have to resort to explicit time stepping methods, or use alternate numerical methods such as pseudospectral methods \cite{Song}.

\subsection{Implementation Concerns} One may question why we have decided to use the Godunov scheme, since for the purposes of implementation, something like the Lax-Friedrichs numerical Hamiltonian may be easier. Indeed, the Lax-Friedrichs Hamiltonian in given by \begin{equation} \label{eq:LFHamiltonian}\hat H_{LF}(u^+_x, u^-_x, u^+_y, u^-_y) = H\left(\frac{u_x^++u_x^-}{2},\frac{u_y^++u_y^-}{2}\right) - \frac{\alpha_1}{2}(u_x^+-u_x^-) - \frac{\alpha_2}{2}(u_y^+-u_y^-), \end{equation} where $\alpha_1, \alpha_2$ are bounds on the derivatives of $H$ with respect to the first or second arguement, respectively \cite{OsherShu}. The Lax-Friedrichs scheme works by adding numerical diffusion, thus in essence changing equation \eqref{eq:HJ} to \begin{equation} \label{eq:LFeq} u_t - \varepsilon \Delta u + H(u_x, u_y) = 0, \end{equation} where $\varepsilon = O(\Delta x, \Delta y).$ In most applications, this is acceptable, since the Lax-Friedrichs Hamiltonian still provides a first order accurate approximation to the Hamilton-Jacobi equation. However, in our case, adding diffusion at level $\varepsilon$ to the Hamilton-Jacobi equation is akin to adding uncertainty in the equation of motion at level $\varepsilon^{1/2}$. For example, in the discretization we will use, the numerical diffusion would be on the order of $0.01$, which would correspond to uncertainty in the equation of motion on the order of $0.1$ m/s. This is a nontrivial level of uncertainty, representing roughly one tenth of the maximum walking velocity. For this reason, minimally diffusive schemes are necessary for our application. A similar roadblock arises when using level set methods, especially when the geometry of the level sets is a crucial aspect of the model like in shock-capturing, image-processing \cite{Fedkiw2003} or recent deforestation models \cite{Deforest}. In these cases, numerical diffusion becomes noticeable when level set velocity is near zero. 

If the Hamiltonian $H$ is relatively simple, the minima and maxima in the Godunov Hamiltonian can be resolved exactly. For example, in the special case that $H(u,v) = h(u^2,v^2)$ and $h$ is monotone in each argument, the Godunov scheme simplifies significantly \cite{OptPathPaper}. In our application, this situation arises when the elevation is flat so that $\nabla E \equiv 0$, and our Hamiltonian \eqref{eq:Hamiltonian} becomes $H(x, \nabla u) = -V_{\text{max}} \|\nabla u \|_2$. In general, our Hamiltonian is much more complicated. In order to implement the Godunov scheme for our Hamiltonian, we must resolve three minima or maxima: the minimum involved in the definition of the Hamiltonian, and the two minima/maxima involved in the scheme itself. We resolve all this minima and maxima discretely, by simply sampling points and choosing the correct one. As long as the error from this descrete optimization remains on the order of $\Delta x$ and $\Delta y$, the scheme will remain accurate to first order. The level of resolution needed for the discrete optimization depends somewhat on the problem, but empirically it appears that the most important facet is the regularity of the elevation data $E$. This makes intuitive sense: for less smooth elevation, the minimization in \eqref{eq:Hamiltonian} taken with respect to the walking direction will require finer resolution to resolve. Likewise, for less smooth elevation, the discontinuities in the derivative of the solution $u$ become more severe. Thus the optimization sets in the Godunov scheme, which have the form $I(u^+_x,u^-_x)$ and $I(u^+_y,u^-_y)$, become larger.  

\section{Simulations \& Results}\label{results} We implemented our model in MATLAB using the numerical schemes (\ref{eq:explicitScheme},~\ref{eq:semiImplicitScheme}) to solve the Hamilton-Jacobi-Bellman equations (\ref{eq:HJB},~\ref{eq:SHJBOurs}), and the forward Euler method to solve the equation of motion \eqref{eq:actualODEStoch}. In the following figures, we will display elevation contours ranging from blue signifying low elevation to yellow signifying high elevation. The starting point $x_0$ will be marked with a green star and the desired end point $x_{\text{end}}$ will be marked with a red star. The lines representing the walking paths will be plotted in colors ranging from green, symbolizing simulations with smaller $\sigma$ values, to red, symbolizing larger $\sigma$ values.  

\subsection{Synthetic Elevation Data} \label{syntheticData} 
\begin{figure}[htbp] 
\centering
\begin{subfigure}[t]{0.31\textwidth}
	\includegraphics[width=\textwidth,trim = 42 24 25 19, clip]{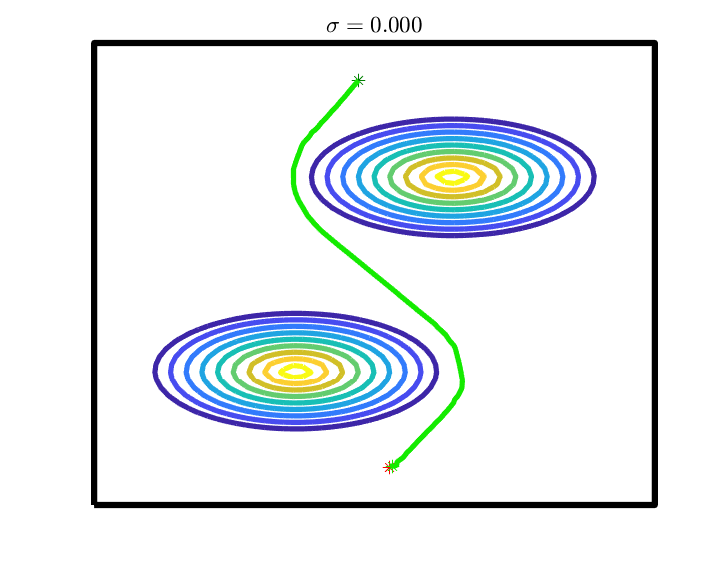}
	\caption{$\sigma = 0.0$}
	\label{fig:twoMountSig0}
\end{subfigure}~
\begin{subfigure}[t]{0.31\textwidth}
	\includegraphics[width=\textwidth,trim = 42 24 25 19, clip]{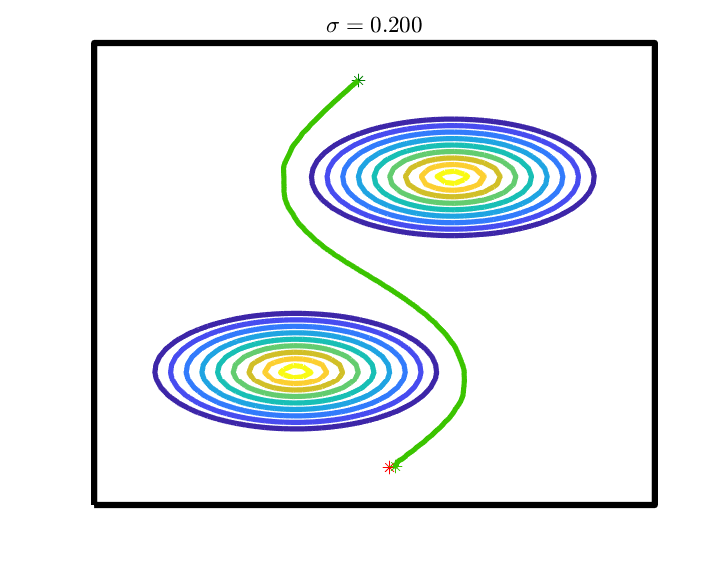}
	\caption{$\sigma = 0.2$}
	\label{fig:twoMountSig2}
\end{subfigure}~
\begin{subfigure}[t]{0.31\textwidth}
	\includegraphics[width=\textwidth,trim = 42 24 25 19, clip]{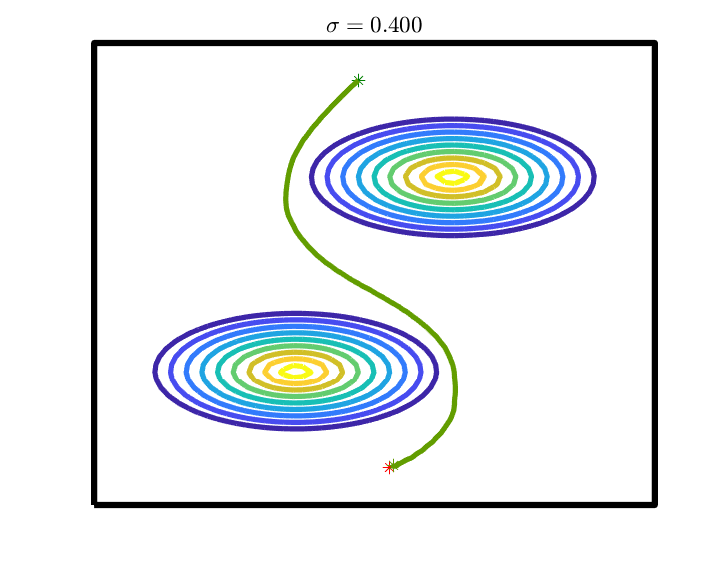}
	\caption{$\sigma = 0.4$}
	\label{fig:twoMountSig4}
\end{subfigure}\\
\centering
\begin{subfigure}[t]{0.31\textwidth}
	\includegraphics[width=\textwidth,trim = 42 24 25 19, clip]{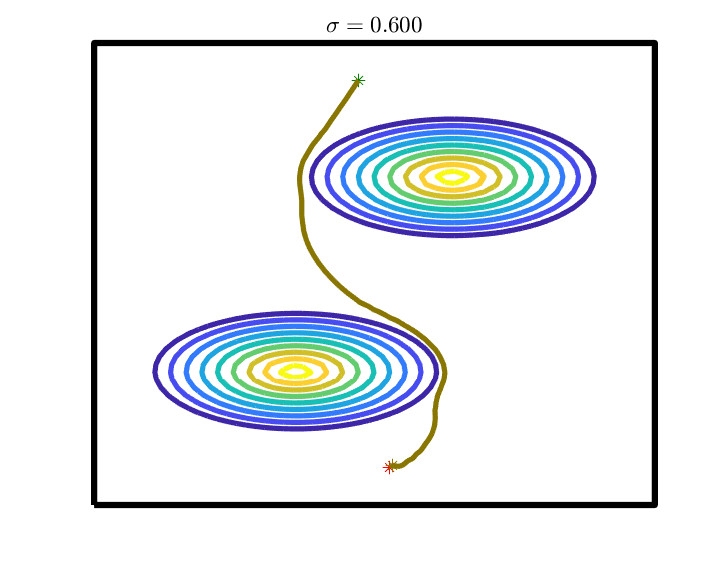}
	\caption{$\sigma = 0.6$}
	\label{fig:twoMountSig6}
\end{subfigure}~
\begin{subfigure}[t]{0.31\textwidth}
	\includegraphics[width=\textwidth,trim = 42 24 25 19, clip]{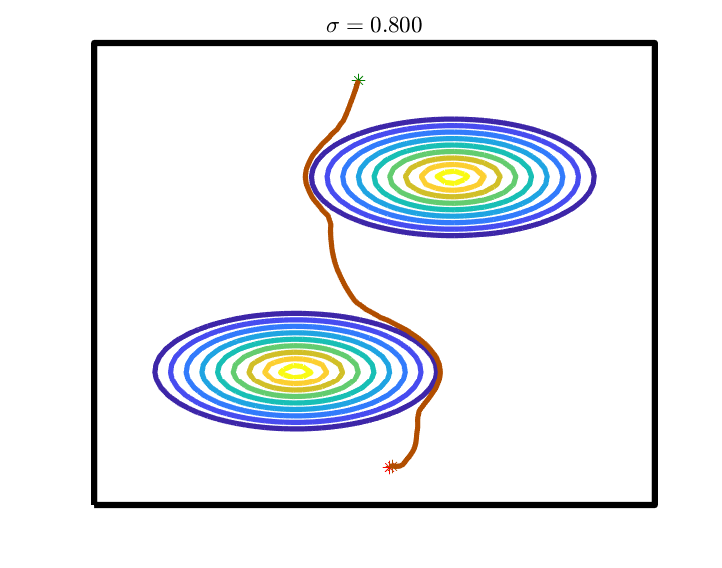}
	\caption{$\sigma = 0.8$}
	\label{fig:twoMountSig8}
\end{subfigure}~
\begin{subfigure}[t]{0.31\textwidth}
	\includegraphics[width=\textwidth,trim = 42 24 25 19, clip]{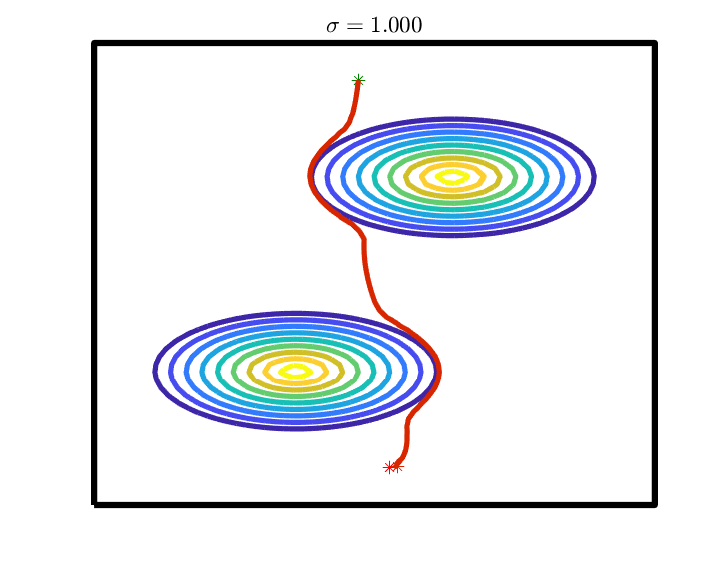}
	\caption{$\sigma = 1.0$}
	\label{fig:twoMountSig10}
\end{subfigure}
\caption{Optimal paths using different $\sigma$ values. End time $T = 3.8$ for each plot.}
\label{fig:twoMountDiffSig}
\end{figure}

We began by testing our model against simple synthetic data. In figure~\ref{fig:twoMountDiffSig}, we have computed optimal paths with several different levels of uncertainty $\sigma$. Referring to table~\ref{fig:twoMethods}, we are using method $(i)$ to compute the optimal paths. That is, we are using the stochastic HJB equation to determine the optimal control values, but then computing the path using the deterministic equation of motion. In this example, the elevation is flat except for two large mountains which lie between the starting point and end point.

In the deterministic case, plotted in figure~\ref{fig:twoMountSig0}, the path curves around the mountains as one would expect: the walking velocity is significantly hampered by the change in elevation, so it is more efficient to avoid those regions. In this case, we see that the optimal path suggested by our algorithm is not particularly sensitive to small changes in $\sigma$. The path in figure~\ref{fig:twoMountSig2} which has $\sigma = 0.2$ looks very similar to that in figure~\ref{fig:twoMountSig0} which has $\sigma = 0.$ However, as $\sigma$ becomes larger, we do see significant changes in the path. The path in figure~\ref{fig:twoMountSig10} where $\sigma = 1$ is significantly different from the deterministic case. Here the uncertainty is on roughly the same order as the walking velocity. With this level of uncertainty, one could imagine walking through a forest in a very dense fog. In planning the path, this algorithm suggests that you walk directly toward your destination and adjust as necessary when obstacles arise. 

\begin{figure}[htbp]
\centering
\begin{subfigure}[t]{0.31\textwidth}
	\includegraphics[width=\textwidth]{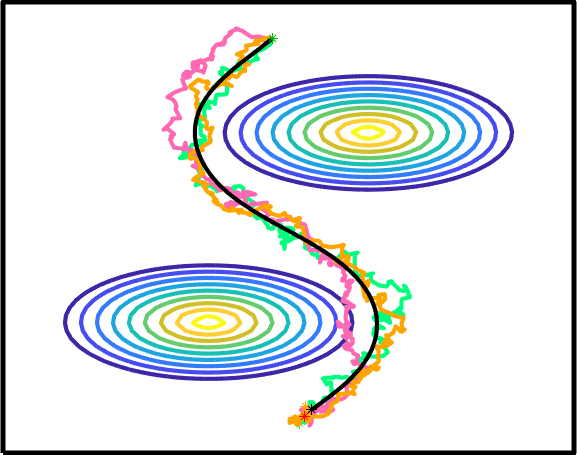}
	\caption{$\sigma = 0.2$}
	\label{fig:stochPathsSig02}
\end{subfigure}~
\begin{subfigure}[t]{0.31\textwidth}
	\includegraphics[width=\textwidth]{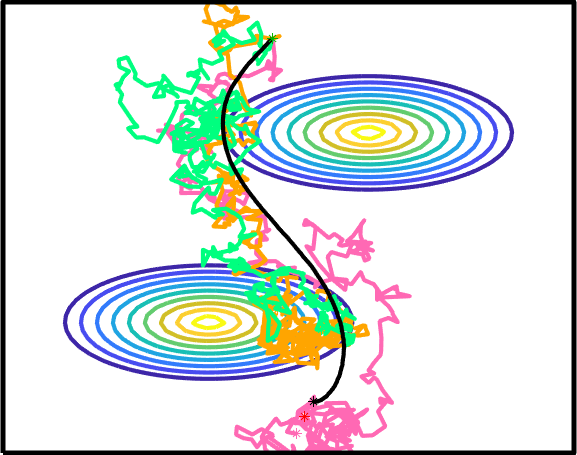}
	\caption{$\sigma = 0.5$}
	\label{fig:stochPathsSig05}
\end{subfigure}~
\begin{subfigure}[t]{0.31\textwidth}
	\includegraphics[width=\textwidth]{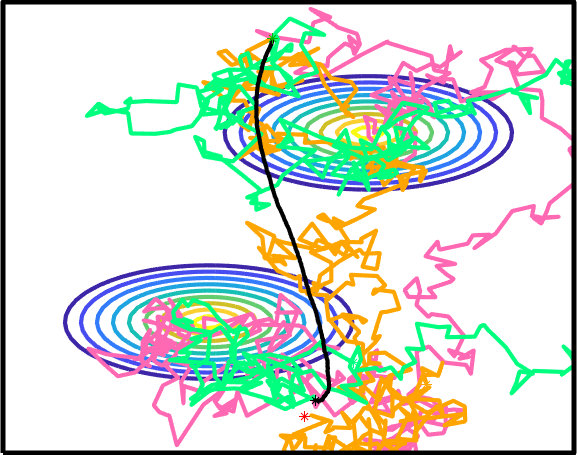}
	\caption{$\sigma = 0.8$}
	\label{fig:stochPathsSig08}
\end{subfigure}\\
\centering
\begin{subfigure}[t]{0.31\textwidth}
	\includegraphics[width=\textwidth]{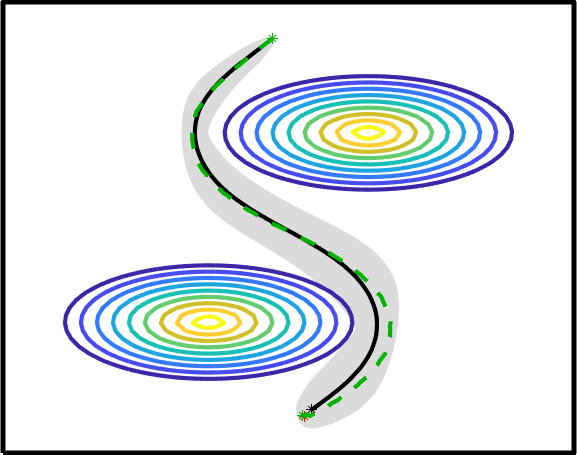}
	\caption{$\sigma = 0.2$}
	\label{fig:varianceSig02}
\end{subfigure}~
\begin{subfigure}[t]{0.31\textwidth}
	\includegraphics[width=\textwidth]{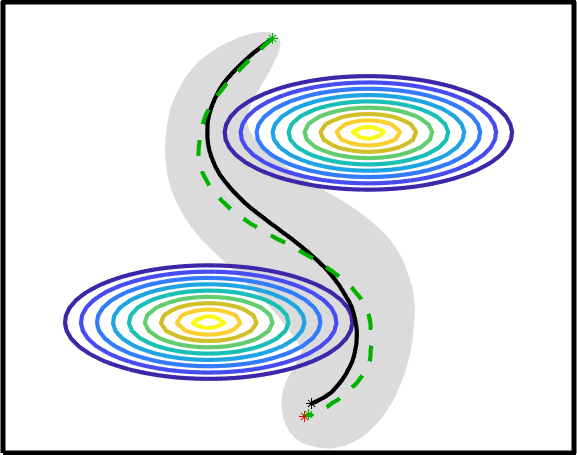}
	\caption{$\sigma = 0.4$}
	\label{fig:varianceSig04}
\end{subfigure}~
\begin{subfigure}[t]{0.31\textwidth}
	\includegraphics[width=\textwidth]{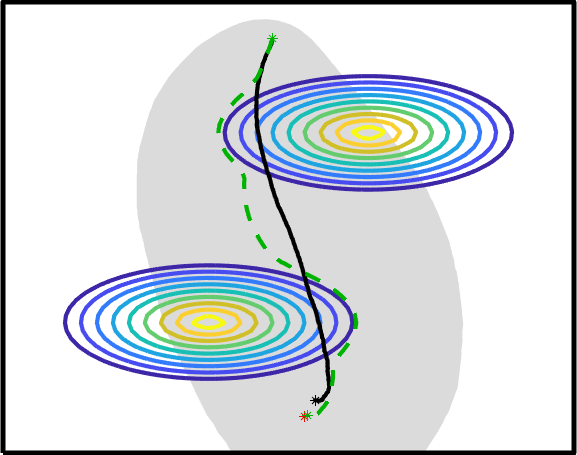}
	\caption{$\sigma = 0.8$}
	\label{fig:varianceSig08}
\end{subfigure}~
\caption{(a) - (c) Average path over 10000 trials (black), and three realizations of the stochastic equation of motion (colored) at different levels of uncertainty $\sigma$. (d) - (f) Average path (black) with standard deviation (grey), and the path computed using method $(i)$ (dotted green).}
\label{fig:twoMountStochPaths}
\end{figure}

Next, we consider method $(ii)$ from table~\ref{fig:twoMethods}; that is, we simulate the stochastic equation of motion many times and compute the average path. Specifically, we simulate the equation $L$ times, resulting in paths $\{\bx_\ell(t)\}_{\ell=1}^{L}$ which are resolved at the same discrete time points, but with some randomness due to the Brownian motion. We then define the average path $\overline \bx(t) = \frac 1 L \sum^L_{\ell=1} \bx_\ell(t)$. We are still using the forward Euler method for the stochastic ODE and since the coefficient in front of the Brownian motion is independent of $X_t$, this corresponds with the Milstein method which exhibits strong and weak convergence at first order \cite{Kloeden1992}. In each of the following results, we have simulated the equation of motion 10000 times. Results are displayed in figures~\ref{fig:stochPathsSig02}-\ref{fig:stochPathsSig08}. The black line represents the average optimal path and the colored lines represent three individual realizations of the stochastic equation of motion. Here as $\sigma$ gets larger, the individual realizations become less meaningful, but the average path is still somewhat smooth and roughly connects the starting point to the ending point. 

We also calculate a form of confidence interval to evaluate how close a single realization is likely to be to the average. To do this, at each point $(\overline x, \overline y)$ along the average path, we calculate the standard deviation $(\delta_x, \delta_y)$ in each of the coordinates. Then at each point, we plot in light grey the ellipse centered at $(\overline x, \overline y)$ with radii $(\delta_x,\delta_y)$ in the $x$ or $y$ direction, respectively. As we travel along the path plotting these ellipses, the grey envelope represents the set of points within one standard deviation of the average path. This is seen in figures~\ref{fig:varianceSig02}-\ref{fig:varianceSig08}. In these plots the average path is plotted as a solid black line. Now we also display the path which was computed using method $(i)$ using a dotted green line. For small $\sigma$, the average path and the determistic path match fairly well. For larger $\sigma$, they begin to diverge, but the walking strategy seems similar: for larger $\sigma$, the average paths take a much more direct approach, cutting corners more closely, or walking directly over the mountains. In each case, the deterministic path stays well within one standard deviation of the average path. Notice that as $\sigma$ gets larger, the standard deviation grows very quickly so that in figure~\ref{fig:varianceSig08} the set of possible paths within one standard deviation of the average is quite large, and it may simply be that method $(i)$ provides a more reasonable solution in this application.

\subsection{Real Elevation Data} \label{realData} 
\begin{figure}[t!]
\centering
\begin{subfigure}[t]{0.46\textwidth}
	\includegraphics[width=\textwidth]{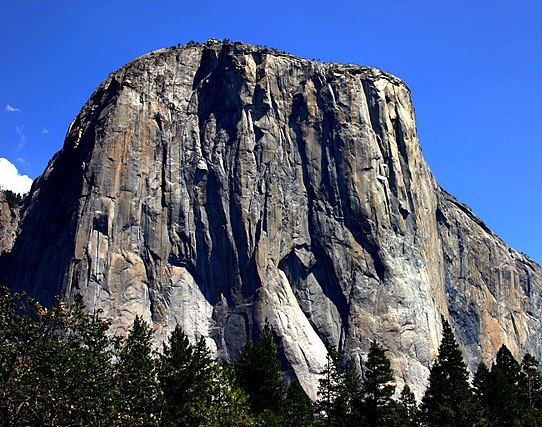}
	\caption{The south-facing cliff face of El Capitan}
	\label{fig:elCapPic}
\end{subfigure}~
\begin{subfigure}[t]{0.46\textwidth}
\includegraphics[width=\textwidth]{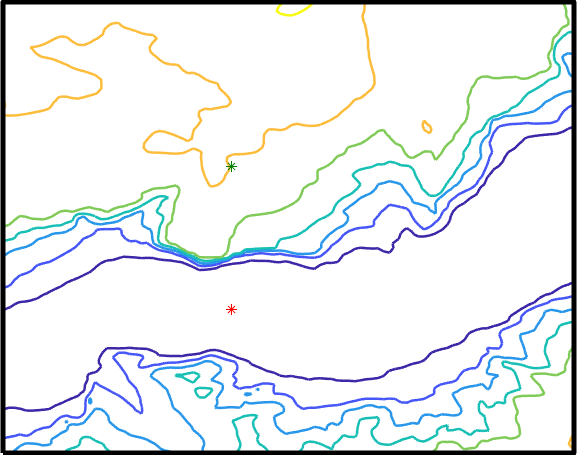}
\caption{The elevation profile of El Capitan.}
\label{fig:elCapElev}
\end{subfigure}
\caption{El Capitan, Yosemite National Park, California\protect\footnotemark} 
\label{fig:elCap}
\end{figure}

Seeing that our model works correctly for simplified elevation data, we tested the model against real elevation in the area surrounding the mountain El Capitan in Yosemite National Park. The elevation profile of El Capitan, along with the starting an ending points, is pictured in figure~\ref{fig:elCap}. Notice that directly in between the starting and ending points, the contour lines lie close together, indicating a sheer cliff face. The starting point is near the summit of the mountain, and the ending point is in the valley to the south of the mountain, so any walking path should choose the gentler grades to the east or west of the cliff face. 

\begin{figure}[t!] 
\centering
\begin{subfigure}[t]{0.31\textwidth}
	\includegraphics[width=\textwidth,trim = 42 24 25 19, clip]{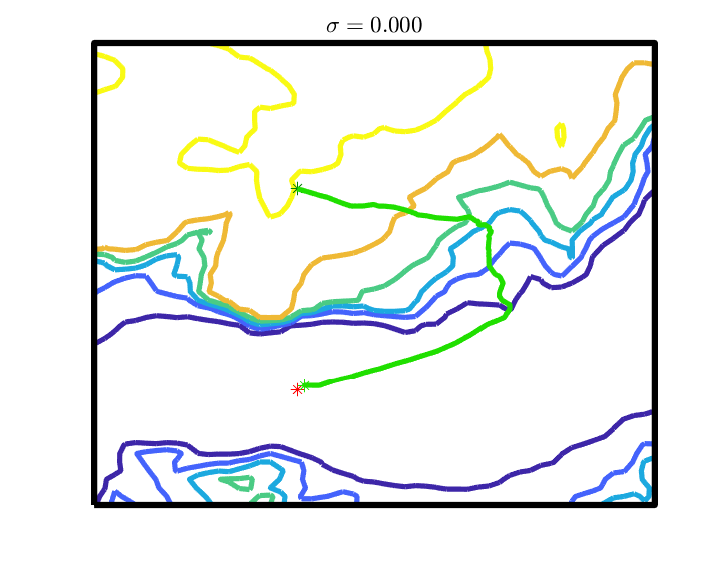}
	\caption{$\sigma = 0.00, \,\, T \approx 6800$}
	\label{fig:elCapSig0}
\end{subfigure}~
\begin{subfigure}[t]{0.31\textwidth}
	\includegraphics[width=\textwidth,trim = 42 24 25 19, clip]{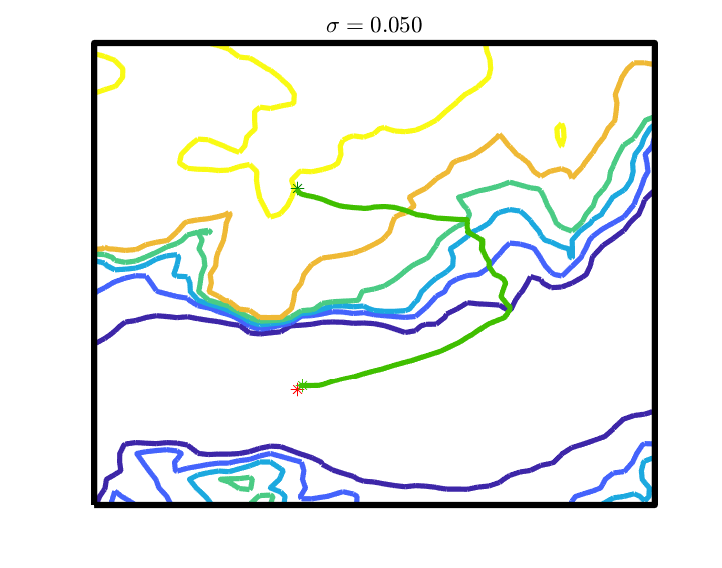}
	\caption{$\sigma = 0.05, \,\, T \approx 6850$}
	\label{fig:elCapSig05}
\end{subfigure}~
\begin{subfigure}[t]{0.31\textwidth}
	\includegraphics[width=\textwidth,trim = 42 24 25 19, clip]{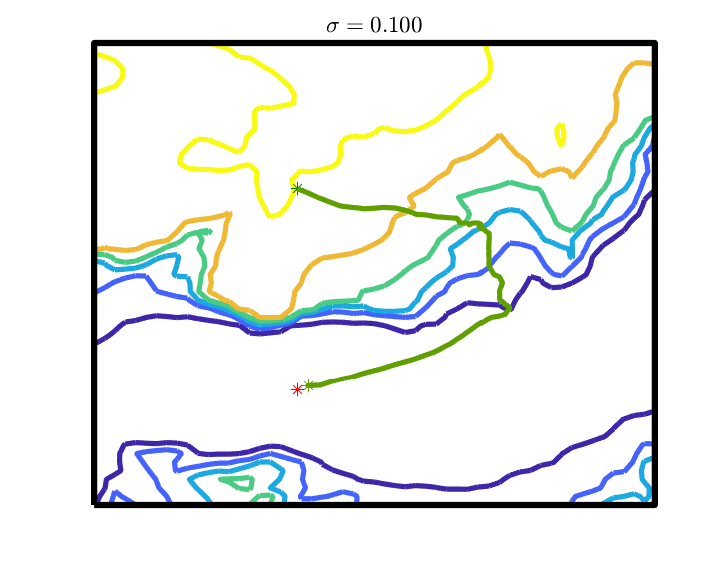}
	\caption{$\sigma = 0.10, \,\, T \approx 7500$}
	\label{fig:elCapSig10}
\end{subfigure}\\
\centering
\begin{subfigure}[t]{0.31\textwidth}
	\includegraphics[width=\textwidth,trim = 42 24 25 19, clip]{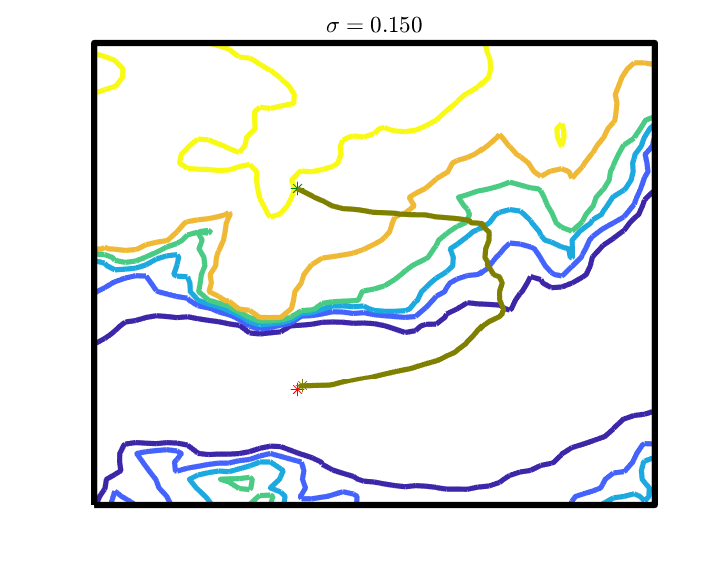}
	\caption{$\sigma = 0.15, \,\, T \approx 7550$}
	\label{fig:elCapSig15}
\end{subfigure}~
\begin{subfigure}[t]{0.31\textwidth}
	\includegraphics[width=\textwidth,trim = 42 24 25 19, clip]{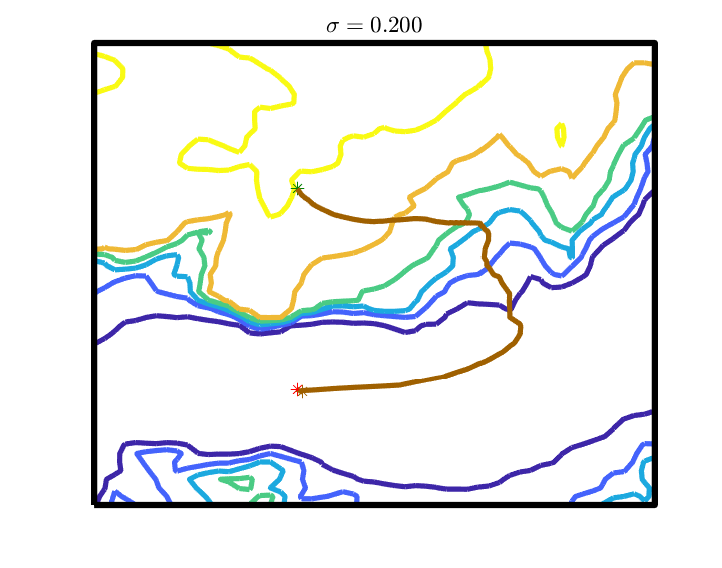}
	\caption{$\sigma = 0.20, \,\, T \approx 8950$}
	\label{fig:elCapSig20}
\end{subfigure}~
\begin{subfigure}[t]{0.31\textwidth}
	\includegraphics[width=\textwidth,trim = 42 24 25 19, clip]{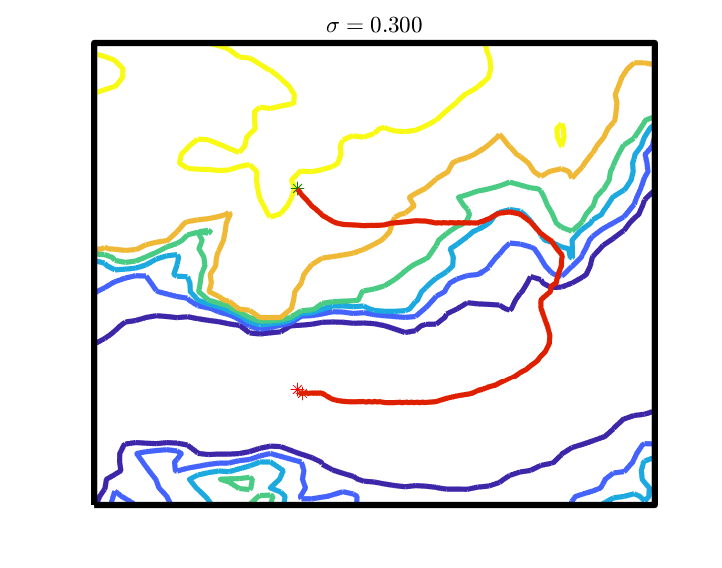}
	\caption{$\sigma = 0.30, \,\, T \approx 15550$}
	\label{fig:elCapSig30}
\end{subfigure}
\caption{Optimal paths descending El Capitan using different $\sigma$ values.}
\label{fig:elCapDiffSig}
\end{figure}

Indeed, this is exactly what we observe, as seen in figure~\ref{fig:elCapDiffSig}. These paths were determined using method $(i)$, the deterministic equation of motion. In these simulations all the scales in the problem are completely genuine. The region displayed in these figures is a rectangle roughly $5$ kilometers east-to-west and $6$ kilometers north-to-south. The starting and ending points are roughly $2$ kilometers apart and the terrain is mountainous, and so several thousand seconds are required to traverse a path connecting the points. Here, the elevation data is much less smooth, and this leads to a greater sensitivity to small changes in $\sigma$. Figure~\ref{fig:elCapSig05} shows that at $\sigma = 0.05$, the optimal path looks largely the same as in the deterministic case, displayed in figure~\ref{fig:elCapSig0}. However, when $\sigma = 0.3$, the optimal path is quite different, as seen in figure~\ref{fig:elCapSig30}.  \footnotetext{Image courtesy of Mike Murphy, uploaded to Wikipedia Commons under Share-Alike license: https://commons.wikimedia.org/wiki/File:Yosemite\_El\_Capitan.jpg} 

One significant note here: at different levels of $\sigma$, there are different optimal terminal times $T$. Recall, the parameter $T$ must be chosen before simulating the model. In the case of the synthetic data in figure~\ref{fig:twoMountDiffSig}, the terminal time $T$ is not particularly sensitive to changes in $\sigma$, since qualitatively the paths are all similar and the amount of time that is ``wasted" by taking a non-optimal path is not significant. In that case, we set $T = 3.8$ and any path with $\sigma \in [0,1]$ had sufficient time to reach the endpoint. This is not the case in figure~\ref{fig:elCapDiffSig}, where small changes in $\sigma$ lead to more significant qualitative changes in the paths. Indeed, the greedy strategy of taking a more direct route and then adjusting as necessary can be very costly in the case of El Capitan where it is very easy to get stuck in regions of severe grades, and be nearly unable to move. In this case, if one is reasonably uncertain about walking velocity as in figure~\ref{fig:elCapSig30}, the algorithm suggests one should allow ample time, and take a safer route which more deliberately avoids regions with large changes in elevation. Because this route is significantly different, it requires a terminal time of roughly $T = 15550$ seconds, as opposed to a terminal time closer to $T=6800$ seconds as in figure~\ref{fig:elCapSig0}. For larger values of $\sigma$ (for example $\sigma> 0.5$), the path will not make it down the mountain even given exorbitantly large terminal time, because it will walk too close to the cliff, become stuck, and have insufficient time to adjust. We say more about the role of the parameter $T$, especially as it pertains to impassable obstacles such as the El Capitan cliff face, in section~\ref{theRoleOfT}. 

\begin{figure}[b!]
\centering
\begin{subfigure}[t]{0.31\textwidth}
	\includegraphics[width=\textwidth]{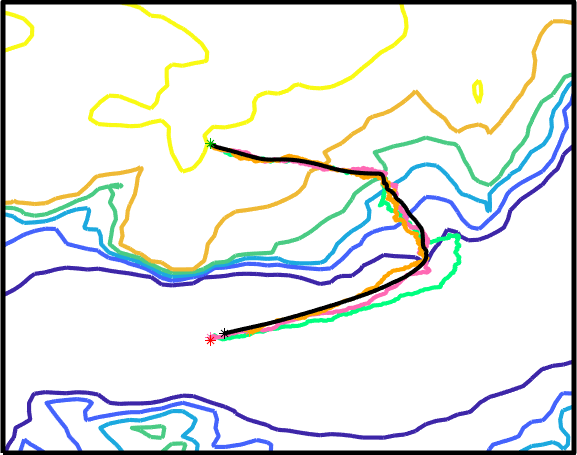}
	\caption{$\sigma = 0.05$}
	\label{fig:stochPathsElCapSig005}
\end{subfigure}~
\begin{subfigure}[t]{0.31\textwidth}
	\includegraphics[width=\textwidth]{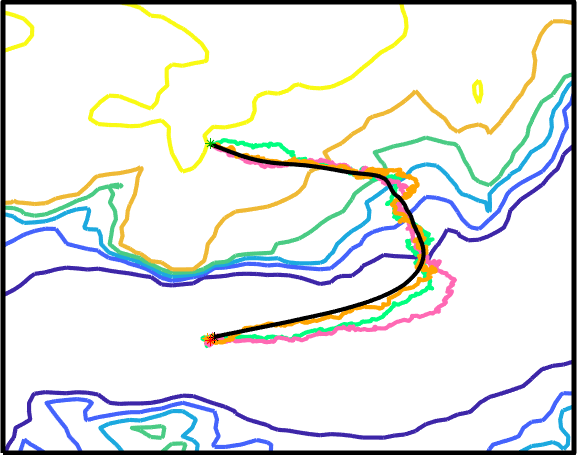}
	\caption{$\sigma = 0.1$}
	\label{fig:stochPathsElCapSig010}
\end{subfigure}~
\begin{subfigure}[t]{0.31\textwidth}
	\includegraphics[width=\textwidth]{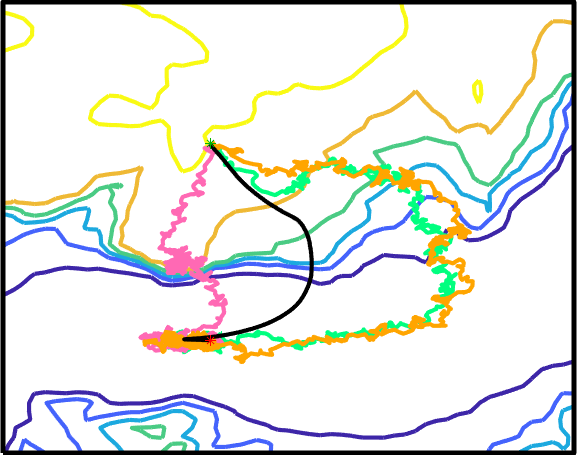}
	\caption{$\sigma = 0.2$}
	\label{fig:stochPathsElCapSig020}
\end{subfigure}\\
\begin{subfigure}[t]{0.31\textwidth}
	\includegraphics[width=\textwidth]{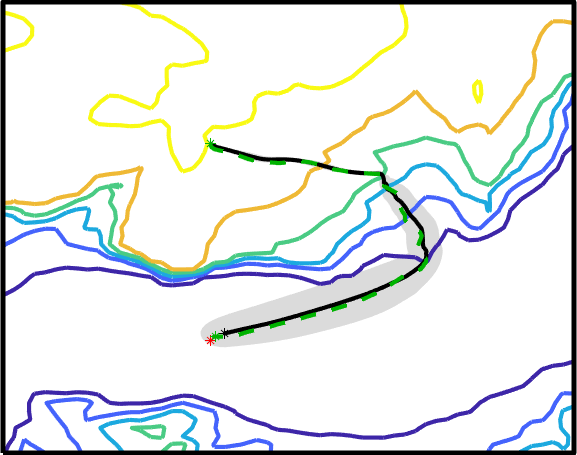}
	\caption{$\sigma = 0.05$}
	\label{fig:varianceElCapSig005}
\end{subfigure}~
\begin{subfigure}[t]{0.31\textwidth}
	\includegraphics[width=\textwidth]{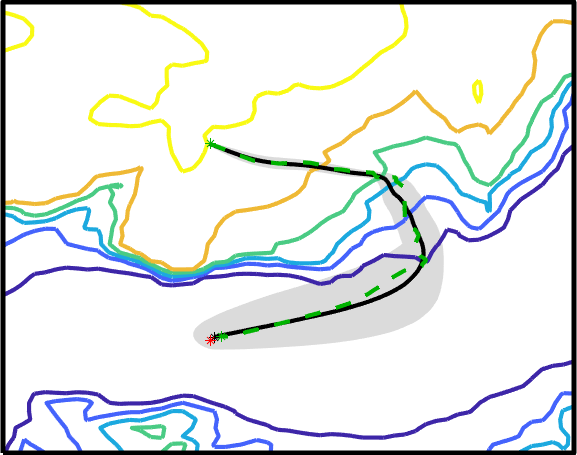}
	\caption{$\sigma = 0.1$}
	\label{fig:varianceElCapSig010}
\end{subfigure}~
\begin{subfigure}[t]{0.31\textwidth}
	\includegraphics[width=\textwidth]{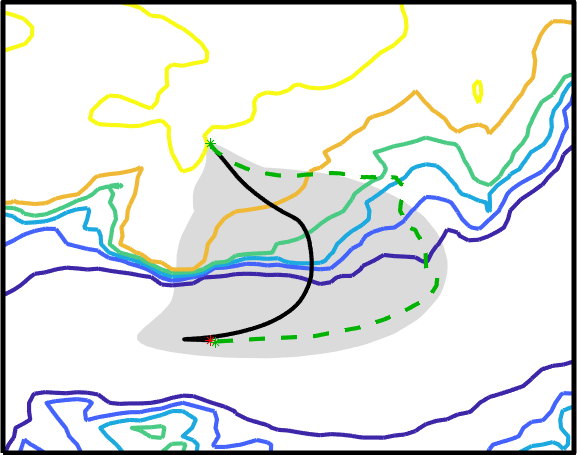}
	\caption{$\sigma = 0.2$}
	\label{fig:varianceElCapSig020}
\end{subfigure}
\caption{(a) - (c) Average path over 10000 trials (black), and three realizations of the stochastic equation of motion (colored) at different levels of uncertainty $\sigma$. (d) - (f) Average path (black) with standard deviation (grey), and the path computed using method $(i)$ (dotted green). }
\label{fig:ElCapStochPaths}
\end{figure}

As in the previous section, we would also like to use method $(ii)$ to construct a path. In figures~\ref{fig:stochPathsElCapSig005}-\ref{fig:stochPathsElCapSig020}, we plot the average path along with three realizations in the case that $\sigma = 0.05, 0.1$ and $0.2$. When $\sigma = 0.05, 0.1$, each of these realizations is fairly close to the average path, and the results are similar to those obtained using method $(i)$. We have also included the region which is one standard deviation away from the average path, as seen in figures~\ref{fig:varianceElCapSig005}-\ref{fig:varianceElCapSig020}. Even when $\sigma$ is very small, the variance in how the paths descend the mountain is fairly large. This is because small perturbations in that region will more qualitatively change the course of path. The results for $\sigma = 0.2$---displayed in figures~\ref{fig:stochPathsElCapSig020},\ref{fig:varianceElCapSig020}---do not seem particularly meaningful. In this case, the uncertainty in the walking velocity is large enough that if the path approaches the large cliff face, the random perturbation can cause the path to move down the cliff. In this region, the walking velocity is approximately zero, and so the random effects are the driving force for the movement. In those simulations, a large enough portion of the paths descended the cliff in this manner, leading to a skewed average path, and an enormously large standard deviation. Similar problems may arise whenever there are regions where the walking velocity is very small. In such cases, it seems that using the deterministic equation of motion with the stochastic control values (as in figure~\ref{fig:elCapDiffSig}) will give a much more meaningful result. 

\subsection{Impassable Obstacles and the Role of the Parameter $T$} \label{theRoleOfT} 
\begin{figure}[b!]
\centering
\begin{subfigure}[t]{0.49\textwidth}
	\includegraphics[width=\textwidth,trim = 42 24 25 19, clip]{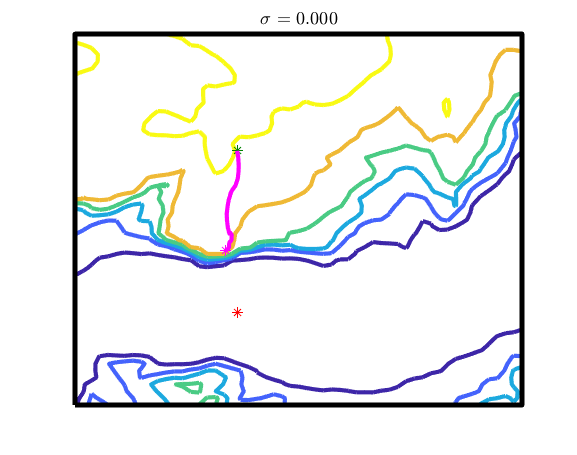}
	\caption{Optimal path given $T = 2000$ seconds.}
	\label{fig:elCapBadPath}
\end{subfigure}~
\begin{subfigure}[t]{0.49\textwidth}
	\includegraphics[width=\textwidth,trim = 42 24 25 19, clip]{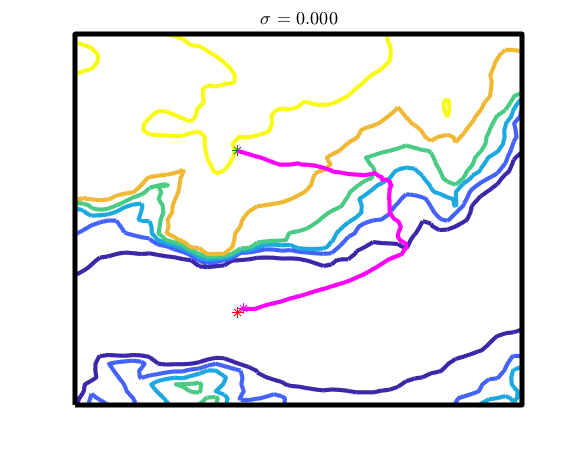}
	\caption{Optimal path given $T = 6800$ seconds.}
	\label{fig:elCapGoodPath}
\end{subfigure}
\caption{Optimal paths in the vicinity of El Capitan with different terminal times.}
\label{fig:elCapDeterministicOpPaths}
\end{figure}

We remarked in section~\ref{HJBDeterministic} that different choices for the terminal time $T$ can lead to qualitative changes in how the path is constructed. We can observe this in the example of El Capitan. In figure~\ref{fig:elCapDeterministicOpPaths}, we used $\sigma = 0$ so that the model is fully deterministic, and we simulated the model with two different terminal times $T$. In figure~\ref{fig:elCapBadPath}, we see that with terminal time $T \approx 2000$ seconds, the path simply walks to the cliff face and stays put. However, given $T \approx 6800$ seconds, the path descends the eastern slope and finds the desired end point as seen in figure~\ref{fig:elCapGoodPath}

\begin{figure}[htbp]
\centering
\begin{subfigure}[t]{0.49\textwidth}
	\includegraphics[width=\textwidth]{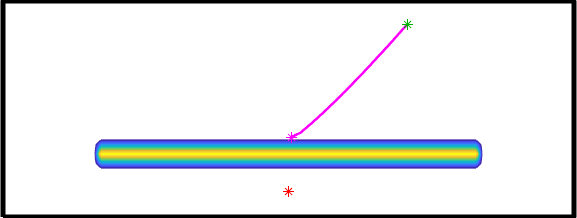}
	\caption{Optimal path given $T = 2$ seconds.}
	\label{fig:largeObs1}
\end{subfigure}~
\begin{subfigure}[t]{0.49\textwidth}
	\includegraphics[width=\textwidth]{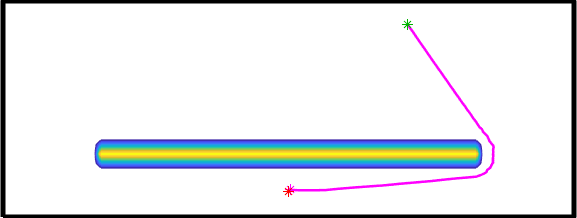}
	\caption{Optimal path given $T = 4.25$ seconds.}
	\label{fig:largeObs2}
\end{subfigure}
\caption{Optimal paths using different end time values. The colored region is a wall.}
\label{fig:largeObs}
\end{figure}

We can recreate this scenario using synthetic elevation data by placing a large wall directly between the starting point and end point as in figure~\ref{fig:largeObs}. The elevation is incredibly steep in the colored region, meaning that any optimal path would surely avoid the wall. In figure~\ref{fig:largeObs2}, where $T = 4.25$, this is exactly the behavior we observe; the path curves around the obstacle. However, in figure~\ref{fig:largeObs1},  where $T = 2$, the path walks toward the obstacle, stopping at the edge because velocity is near zero there. 

Recall, our model constructs the path which ends as close (in Euclidean distance) to the desired end point as possible in the time allotted. When $T = 2$ in figure~\ref{fig:largeObs1}, there is not enough time to walk around the wall and instead, to get as close to the desired end point as possible, the path walks directly toward the wall. When situations like this arise, there is some critical amount of time $T^* > 0$ such that, given $T > T^*$, the path will walk around the obstacle, but given $T < T^*$, the path will walk toward the obstacle because it will not be able to get close enough to the end point if it walks around. 

\begin{figure}[t!]
\centering
\begin{subfigure}[t]{0.49\textwidth}
	\includegraphics[width=\textwidth]{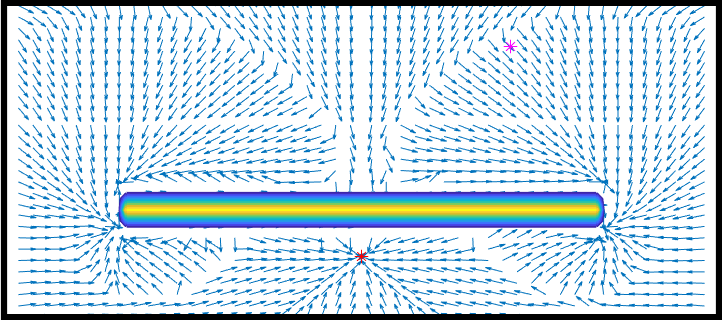}
	\caption{$t = 0.0$}
	\label{fig:sDiscont1}
\end{subfigure}~
\begin{subfigure}[t]{0.49\textwidth}
	\includegraphics[width=\textwidth]{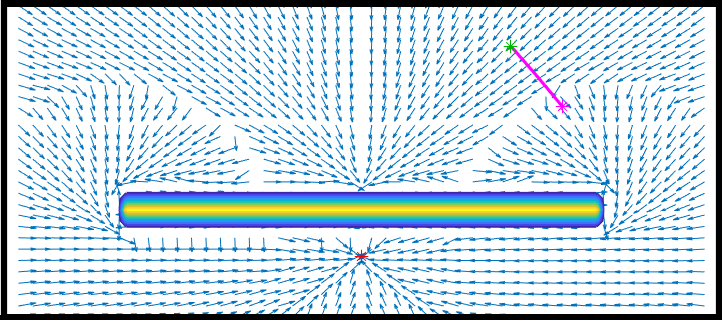}
	\caption{$t = 0.7$}
	\label{fig:sDiscont2}
\end{subfigure} \\
\centering
\begin{subfigure}[t]{0.49\textwidth}
	\includegraphics[width=\textwidth]{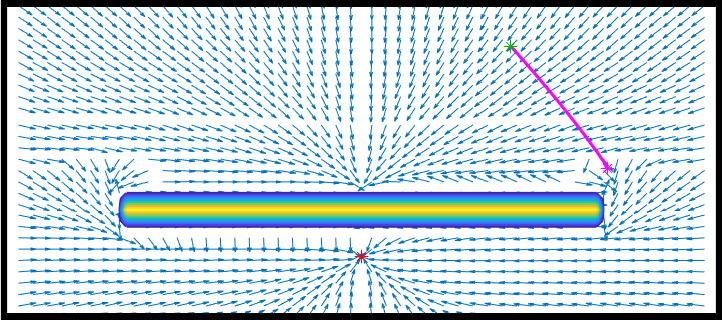}
	\caption{$t = 1.4$}
	\label{fig:sDiscont3}
\end{subfigure}~
\begin{subfigure}[t]{0.49\textwidth}
	\includegraphics[width=\textwidth]{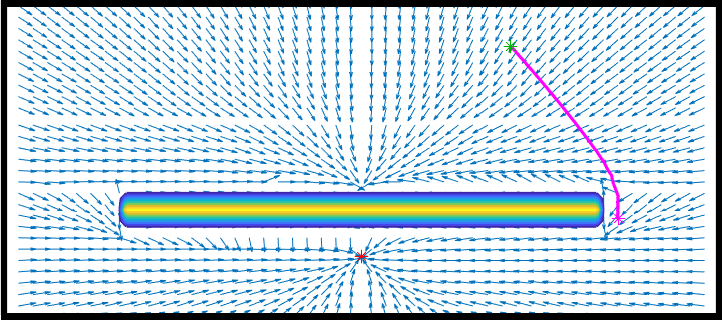}
	\caption{$t = 2.1$}
	\label{fig:sDiscont4}
\end{subfigure}
\caption{The discontinuity in the control value $\B s^*(x,t)$ propagates as time increases}
\label{fig:sDiscont}
\end{figure}

We can see this more explicitly if we plot the actual control values $\B s^*(x,t)$ as well, as is done in figure~\ref{fig:sDiscont}. In this example, the critical time is roughly $T^* = 3.4$, so we have plotted the path created by the algorithm with final time of $T = 3.5$, but we have plotted the path at times $t = 0,0.7,1.4, 2.1$. The arrows in the pictures are the values of $\B s^*(x,t)$. Notice that in figure~\ref{fig:sDiscont1}, there appears to be a discontinuity in the optimal control value. The deciding factor for whether the path will walk around the wall or walk toward the wall is where the starting point lies relative to this discontinuity. As time advances, the discontinuity in $\B s^*(x,t)$ propagates, and since the starting point lies below the discontinuity, the path follows the arrows and walks around the obstacle. In the case when $T = 2$, the starting point is above the discontinuity, and thus the path walks toward the obstacle, rather than around it. 

Discontinuities in $\B s^*(x,t)$ are to be expected and relate to non-uniqueness of the optimal path. If $x_0$ lay directly on the discontinuity in $\B s^*(x,0)$, then either walking around the obstacle or toward it would be equally optimal, since both would result in a path that ends the same distance from the desired end point. Mathematically, one reason to expect discontinuities in $\B s^*(x,t)$ is because $\B s^*(x,t)$ is closely related to the  gradient $\nabla u(x,t)$ of the solution to the HJB equation. Indeed, in the case of isotropic motion (that is, when the equation of motion has the form $\dot {\B x} = f(\B x)\B s$), we will have $\B s^*(x,t) = -\nabla u(x,t) / \| \nabla u (x,t)\|_2$. When the motion is anisotropic, as is the case in our model, the relationship between $\B s^*(x,t)$ and $\nabla u(x,t)$ is no longer so explicit, but we can still anticipate that discontinuities in $\nabla u(x,t)$ will give rise to discontinuities in $\B s^*(x,t)$, and as stated above, we expect the solutions of Hamilton-Jacobi equations to have discontinuities in their derivatives \cite{Bardi1997,CrandallLions,osher2003level}. 

\begin{figure}[t!]
\centering
\begin{subfigure}[t]{0.46\textwidth}
\includegraphics[width=\textwidth]{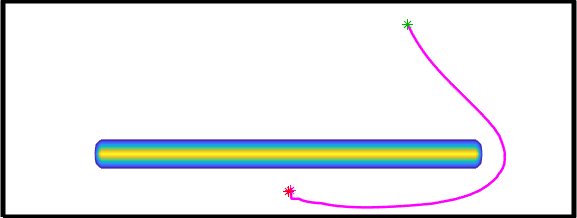}
\caption{Method $(i)$, $\sigma = 0.3$, $T = 4.75$.}
\label{fig:wallStoch1}
\end{subfigure}~
\begin{subfigure}[t]{0.46\textwidth}
	\includegraphics[width=\textwidth]{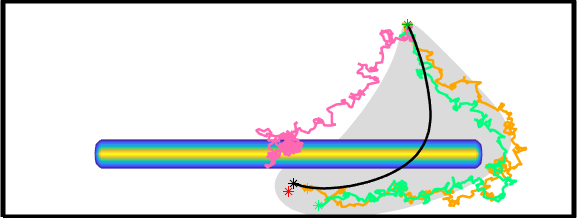}
	\caption{Method $(ii)$, average path, standard deviation and three realizations.}
	\label{fig:wallStoch2}
\end{subfigure}
\caption{Optimal paths around the wall with uncertainty.} 
\label{fig:wallStoch}
\end{figure}

When we add uncertainty, the path planning strategy becomes more greedy, walking directly toward the end point and adjusting to avoid obstacles as is seen in figure~\ref{fig:twoMountDiffSig}. When there is a wall, this strategy is costly because if one walks toward the wall, there may be insufficient time to adjust the route, and thus the critical time $T^*$ required to walk around the wall increases rapidly. This is why the large increase in $T$ is necessary in the example of El Capitan in figure~\ref{fig:elCapDiffSig}. We observe the same behavior in this synthetic example with the wall, though the increase in $T$ as not as pronounced as in the case of El Capitan. In figure~\ref{fig:wallStoch1}, we set $\sigma = 0.3$ and notice that to wrap around the wall and reach the end point, the optimal path computed using method $(i)$ requires an end time of $T = 4.75$ rather than $T = 4.25$ in the deterministic case.  In figure~\ref{fig:wallStoch2}, we use method $(ii)$, computing the average path over 10000 trials, and the path cuts through the wall since enough individual paths were pushed off course due to the random perturbations, as is seen with the pink path in the figure. As in figure~\ref{fig:varianceElCapSig020}, an individual could not realistically traverse the average path, since the wall is impassable. Thus, it seems that method $(i)$ gives to a more practical result. 

The dependence of the model on the parameter $T$ is a major qualitative difference between this optimal path planning model and the model presented in \cite{OptPathPaper}. That model deals only with the fully deterministic case, and uses a level set formulation wherein level sets representing optimal travel evolve outward from the starting point, and the terminal time $T$ is defined as the time required for the level sets to envelop the end point. However, when we add uncertainty to the model, we introduce diffusion in the HJB equation and lose the level set intepretation of the equation. Thus while the model in \cite{OptPathPaper} has the advantage of not depending on $T$, this model is more generally applicable.   

\subsection{Convergence to the Deterministic Path as $\sigma \searrow 0$}\label{convergence} 
\begin{figure}[b!]
\centering
\begin{subfigure}[t]{0.46\textwidth}
\includegraphics[width=\textwidth]{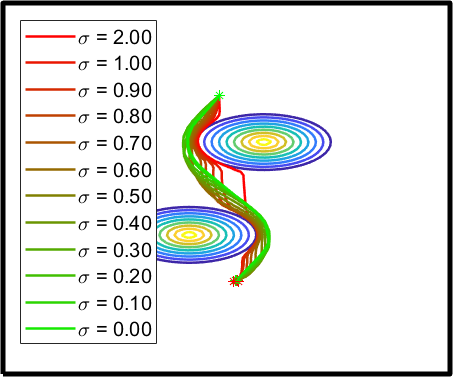}
\caption{Convergence of paths with synthetic elevation data.}
\label{fig:toyDataConvergence}
\end{subfigure}~
\begin{subfigure}[t]{0.46\textwidth}
	\includegraphics[width=\textwidth]{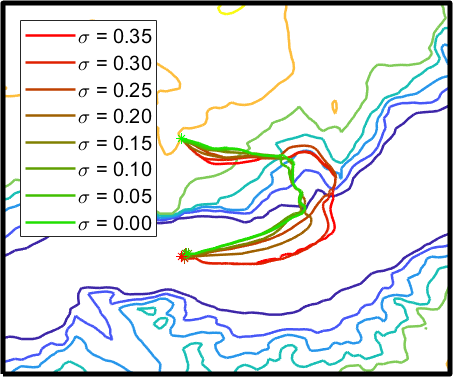}
	\caption{Convergence of paths at El Capitan.}
	\label{fig:elCapConvergence}
\end{subfigure}
\caption{As $\sigma \searrow 0$, the stochastic optimal path converges back to the deterministic optimal path.} 
\label{fig:convergence}
\end{figure}

As stated in section~\ref{ourmodel}, given mild regularity conditions on our Hamiltonian, the solution to the stochastic HJB equation \eqref{eq:SHJBOurs} will converge to the viscosity solution to the ordinary HJB equation \eqref{eq:HJB} as $\sigma \searrow 0$ \cite{CrandallLions}. We can see this empirically, not by observing the solution itself, but by examining the optimal path suggested by our algorithm at different levels of $\sigma$. This is shown in figure~\ref{fig:convergence}. Here we have plotted many paths on the same figure, each computed using method $(i)$ with a different $\sigma$ value. As before, paths plotted in green were computed using smaller $\sigma$ values, and those in red were computed using larger $\sigma$ values. In figure~\ref{fig:toyDataConvergence}, we see a very clear color gradient: the red paths computed with larger $\sigma$ clearly tend toward the green path as $\sigma$ decreases to zero. In figure~\ref{fig:elCapConvergence}, this is less obvious, especially since, for larger $\sigma$, the path takes a qualitatively different route. However, we do see that for smaller $\sigma$ (greener paths), there is a tendency toward the deterministic optimal path.

\section{Conclusions}\label{conclusions} Path planning algorithms have wide-reaching applications in self-driving vehicles, reach-avoid games, pedestrian flow modeling and many other areas. Many previous models for path planning are completely deterministic, while in reality stochastic effects may be present and can significantly alter the motion along the path.   

In this paper we develop a method for optimal path planning of human walking paths in mountainous terrain using a control theoretic approach and a Hamilton-Jacobi-Bellman (HJB) equation and allowing for uncertainty in the controlled equation of motion. The walking speed in our model depends on local slope in the direction of travel, giving rise to an anisotropic control problem. In the HJB equation, the uncertainty presents itself in the form of diffusion, leading to a viscous Hamilton-Jacobi-type equation. We suggest numerical methods for solving these equations, opting for a semi-implicit numerical scheme with a minimally diffusive numerical Hamiltonian, since any spurious numerical diffusion could be intepreted as nontrivial amounts of uncertainty in the equation of motion. After solving the HJB equation numerically, we suggest two methods for resolving the optimal path. First, we use the optimal control values resolved via the stochastic HJB equation, but simulate a deterministic equation of motion. This could simulate a person walking through a dark room or a dense forest, wherein they are cognizant of some uncertainty as they are planning the route, but do not feel random perturbations in velocity as they walk along a path. Second, we integrate the stochastic differential equation many times and arrive at a single path by averaging the results. This could model scenarios such as underwater unmanned vehicles, wherein the traveler actually feels the stochastic effects on the travel velocity. 

We test our algorithm, including both methods for resolving the path, with synthetic elevation data first, and then with real elevation data in the area surrounding El Capitan. We compare these two notions of optimal path, and conclude that in the case of real elevation data or impassable barriers, the first notion gives a more meaningful result. We also observe that in these cases, there will be discontinuities in the optimal control parameter and, especially in the presence of large walls, the position of these discontinuities can determine the walking strategy. Finally, we simulate the model at different levels of uncertainty in the equation of motion and observe that as uncertainty tends to zero, the optimal path path suggested by the model converges back to the deterministic optimal path.

\section{Acknowledgements} Financial support for this work was furnished by the NSF (DMS-1737770). This work was also funded by an academic grant from the National Geospatial-Intelligence Agency (Award No. \#HM0210-14-1-0003, Project Title: Sparsity models for spatiotemporal analysis and modeling of human activity and social networks in a geographic context). Approved for public release, 20-528.

The elevation data for Yosemite National Park was downloaded from the United States Geological Survey national map viewer \cite{NationalMap}. The data was processed and re-formatted using QGIS, an open source geological data processing tool \cite{QGIS}. The data was imported into MATLAB using TopoToolbox \cite{topotoolbox1,topotoolbox2}.

\bibliographystyle{abbrvnat}
\bibliography{bibliography}
\end{document}